\documentclass{article}
\usepackage[utf8]{inputenc}
\usepackage[english]{babel}
\usepackage{csquotes}
\usepackage{biblatex}
\usepackage{hyperref}

\usepackage{ dsfont }
\usepackage{amssymb}
\usepackage{amsmath}        
\usepackage{amsfonts}       
\usepackage{amsthm}         
\usepackage{bbding}         
\usepackage{bm}             
\usepackage{graphicx}       
\usepackage{fancyvrb}       
\usepackage{indentfirst}    
\usepackage{icomma}         
\usepackage{dcolumn}        
\usepackage{booktabs}       
\usepackage{paralist}       
\usepackage{xcolor}         
\usepackage{ dsfont }

\input xy                         
\xyoption{all}
\usepackage{amsmath,amscd}
\usepackage{textcomp}
\theoremstyle{plain}
\newtheorem{veta}{Věta}
\newtheorem{Thm}[veta]{Theorem}
\newtheorem{Prop}[veta]{Proposition}
\newtheorem{Ex}[veta]{Example}

\newtheorem{Lemma}[veta]{Lemma}
\newtheorem{Alg}[veta]{Algorithm}

\theoremstyle{plain}
\newtheorem{Def}[veta]{Definition}

\theoremstyle{remark}

\newtheorem{remark}[veta]{Remark}

\theoremstyle{plain}

\newcommand{\R}{\mathbb{R}}

\newcommand{\N}{\mathbb{N}}

\newenvironment{dukaz}{
  \par\smallskip\noindent
  \textit{Proof}.
}{

\rightline{$\qedsymbol$}
}

\usepackage{blindtext}
\title{Generalised Gabriel-Roiter measure and thin representations}
\author{Dominik Krasula \footnote{This work is a part of project SVV-2023-260721}
}

\begin{document}
 \maketitle 
\begin{abstract}
For Dynkin and Euclidean quivers, it is shown that Gabriel-Roiter measures of thin representations equal the induced chain length functions on the corresponding system of subquivers. This allows a combinatorial procedure to find GR filtrations of thin representations, showing that GR measures of thin representations are field-independent. It is proved that an indecomposable filtration of a thin representation is a GR filtration for a suitable choice of a length function on the category of finite-dimensional representations. 
\end{abstract}

\smallskip


\smallskip

\noindent \textbf{Author:} RNDr. Dominik Krasula

Charles University, Faculty of Mathematics and Physics, Prague, Czech Republic

krasula@karlin.mff.cuni.cz

ORCID: 0000-0002-1021-7364

\medskip

Krasula, D. (2025). Generalised Gabriel-Roiter measure and thin representations. \textit{J. Math.} 663. 

DOI: https://doi.org/10.1016/j.jalgebra.2024.09.017.

\medskip 

\noindent This research was supported by the grant GA ČR 23-05148S from the Czech
Science Foundation.

\smallskip

\section{Introduction}

In 1968, Roiter published proof of the first Brauer-Thrall conjecture in [18]. He showed that a finite-dimensional algebra $A$ over an arbitrary field has an infinite number of pairwise nonisomorphic indecomposable representations (i.e., $A$ is representation-infinite) if and only if it has indecomposable representations of arbitrarily high composition length (i.e., $A$ is of unbounded representation type).  As noted by Roiter, the result holds for any Artinian ring. 

The idea behind Roiter's proof was to show that for each finite-dimensional algebra $A$ of bounded representation type, there exists a function $f$ assigning natural numbers to indecomposable modules of finite length, refining the composition length. Moreover, the value on the image of an epimorphism is bound by the values on direct summands of the domain. The existence of the function is shown by creating a collection of full subcategories of the category of indecomposable representation.

Gabriel, in his 1972 report, [11],  defines \textit{the Roiter measure} for abelian length categories, giving a combinatorial interpretation
of the induction scheme used in Roiter's proof.  He formulated three \textit{basic properties} and proved the\textit{ main property}, implying that a measure of a decomposable object is equal to the maximal measure of its indecomposable summands. All of Gabriel's results about the Roiter measure hold for a general abelian length category. Still, in his article, he works under the assumption that there is an upper bound for lengths of indecomposable objects.  In [16], Ringel speculates that this restriction to the representation-finite algebras explains why the measure was largely ignored in the 20th century. Indeed, in the last two decades, the literature about the measure usually discusses representation-infinite algebras. A notable exception is [6].

The measure is now known as the \textit{Gabriel-Roiter measure},  or the \textit{GR measure} for short, first coined by Ringel in his 2005 article [14]. He also dualised the notion, defining the \textit{Gabriel-Roiter comeasure}.  Ringel also extends Gabriel's definition to arbitrary modules. A sequence of indecomposable modules witnessing the GR measure of a given indecomposable module is called a \textit{Gabriel-Roiter filtrations}. The inclusion of two subsequent modules in such a sequence is called a \textit{GR inclusion}.

In [14], the partition of a category $ind\text-A$ of indecomposable representations of artin algebras was obtained, providing two new proofs of the first Brauer-Thrall conjecture using the Gabriel-Roiter measure. It should be noted that the Gabriel-Roiter measure does not satisfy the properties of the function $f$ from Roiter's proof. The Gabriel-Rotier comeasure does.  The GR measure was recently used to prove the Brauer-Thrall conjecture for complete Cohen-Macaulay local rings in [4].

\medskip 

In their seminal article [3], Auslander and Smal\o ~say that Roiter's proof in [18] and Gabriel's interpretation in [11] were \textit{the original impetus} for the work that resulted in their definition of preprojective and preinjective modules over a general artin algebra. A preprojective partition of a category can be seen as a generalisation of Roiter's result. After Ringel revived interest in the GR measure, the interplay of these two theories was studied in depth. Ringel's partition of $ind\text-A$ was compared with the partition to preprojective, regular and preinjective components in [8] and further investigated in [19].

Bo Chen's work on hereditary algebras in [6] suggests that GR-inclusions are of independent interest as a class of inclusion with particular properties. The connection between GR inclusions and irreducible morphisms was described in [17].   GR inclusions and the corresponding factor modules called \textit{GR factors} appear repeatedly in the work of Bo Chen; see, for example, [7] and [10].

\medskip 

Despite the historical connection with artin algebras and AR theory, the idea behind the Gabriel-Roiter measure is purely combinatorial. Krause formalised this in his 2007 article [12]. He defined an analogue of the GR measure for any poset with a \textit{length function}. He also observed that the \textit{basic properties} formulated by Gabriel can be used to axiomatise the Gabriel-Roiter measure. 

Using Krause's approach, one can define different GR measures based on different length functions. Most of the classical results on  GR inclusions, including Ringel's partition of a module category, hold for a generalised GR measure with little or no adjustment to proofs. On the other hand, these results are usually untrue for a general inclusion of two indecomposable modules containing no indecomposable module between them.

Following [19], a sequence of indecomposable submodules of an indecomposable module that cannot be further refined is called \textit{indecomposable filtration}. Theorem \ref{mainThm} shows that if $Q$ is a quiver whose underlying graph is either a tree or $\tilde{A}_n$, then for any indecomposable filtration $\mathcal{F}$ of a thin representation of $Q$, there exists a length function such that $\mathcal{F}$ is a GR filtration for the corresponding GR measure. The necessity of the assumption that the representation is thin is discussed in Example \ref{D4}. 

The key idea in the proof is the realisation of thin representations as subquivers and their lengths as weights of the corresponding quivers. Using this translation, the main theorem is then a straightforward corollary of a combinatorial Theorem \ref{filtration} about length functions on a system of finite subsets. This more general statement allows for more transparent proof and can be of independent interest.

Working with subquivers also gives a simple greedy algorithm for calculating the GR measure of a thin representation over a general field. A natural setting of our methods is quivers whose underlying graph is a tree. However, an extension is made to the case $\tilde{A}_n$.  As the algorithm works only with subquivers, it follows that the GR measure of a given thin representation and corresponding GR filtrations are field-independent.   Analogous results have been proved using Hall polynomials. In particular, it holds for representations of Dynkin quivers, [15], and indecomposable preprojective or quasi-simple homogenous representations of tame quivers; see [20]

So far,  known methods for calculating the GR measure are based on the Auslander-Reiten theory and only work over algebraically closed fields. In this setting, GR measures of indecomposable regular representations of Euclidean quivers can be calculated using results from [8, Section 4]. The case $\tilde{A}_n$ is well-understood; see [9]. In [19], using that a path algebra of a quiver of type $\tilde{A}_n$ is a string algebra,  a combinatorial algorithm is presented to calculate GR measures (w. respect to composition length) of preprojective and non-homogenous regular modules.  If $A$ is an indecomposable representation-finite algebra, the GR measure of all indecomposable modules can be calculated directly from the (finite and connected) Auslander-Reiten quiver of $A$ using  [17, Theorem A]; see Example \ref{D4}.  

\medskip 

The paper is organised as follows. Section \ref{ChainSec} recalls Krause's definition of induced chain length function from [12] and discusses filtrations in a system of finite subsets, Theorem \ref{filtration}. Brief Section \ref{Abel} defines the GR measure as a particular case of an induced chain length function and slightly reformulates Krause's definition into a computational-friendly form.  Section \ref{ThinSec} then describes thin representations in terms of their support quivers.  The main theorem is proved in Section \ref{mainSec} followed by a discussion of the applications and limits of the results.

\section{Chain length functions} \label{ChainSec}

This section recalls Krause's definition of \textit{induced chain length} function from [12], generalising the GR measure for general posets. An analogue of GR filtration is defined. Further examples and properties can be found in [12]. The section ends with a theorem about filtrations on systems of subsets.  

\smallskip 

 Let $(S,\leq)$ be a poset. A finite subset
        $X\subseteq S$ is a \textit{finite chain in S} if the partial order $\leq$  restricted to $X$ is a total order. We write $Ch(S)$ for the set of all finite chains in S.
        
        Let $x,x'\in S$ such that $x<x'$ and for  any element $y$ if $x\leq y \leq x'$ then either $x=y$ or $x'=y$. Then $x$ is called \textit{a direct predecessor of} $x'$ and $x'$ is a \textit{direct successor} of $x$. 

        Let $n\in \N$ and let  $X=\{x_i \mid i\leq n\}\in Ch(S)$ be a finite chain such that $x_i < x_j$ iff $i < j$. We say that $X$ is a \textit{filtration of $x_n$} if $x_1$ is a minimal element of $(S,\leq)$ and for $i<n$ the element $x_i$ is a direct predecessor of $x_{i+1}$.    We  denote filtrations in the following form $\mathcal{X}=(x_1 <  x_2 < \dots <  x_{n-1} <  x_n)$. For $i\leq n$ we denote by $\mathcal{X}^i$ the induced filtration $(x_1 <  x_2 < \dots <  x_{i})$ of $x_i$.

        For a finite chain $X$,  $\max(X)$ ($\min(X)$) denotes the maximum (minimum) of $X$. We use the convention that $\max(\emptyset)  < x  < \min (\emptyset)$  for any $x\in S$

       On $Ch(S)$ we consider the \emph{lexicographical order} defined by
\[X\leq_L Y \Leftrightarrow \min(Y\setminus X) \leq \min(X\setminus Y).\]

Any chain $X \in  Ch(S)$ can be interpreted as a string of 0s and 1s which is indexed by the elements in $S$. The usual lexicographic order on such strings coincides with the lexicographic order on $Ch(S)$. 

\begin{Def}\label{ldef}
       Let $(S,\leq)$ and $(T,\lesssim)$ be two posets.

    Then map  $l\colon S\to T$ is a \emph{length function on $(S,\leq)$} if for any $x,y\in S$:

    (L1) $x< y$ implies $l(x)<(y)$.  

    (L2) $l(x)\lesssim l(y)$ or $l(y)\lesssim l(x)$. 

    (L3) $\{l(x')\mid x' \in S;~ x'\leq y\} $ is finite.   
\end{Def}
If $(T,\lesssim)$ is isomorphic to $(\N,\leq )$, then $(2)$ and $(3)$ always hold. 
\begin{Ex} \label{evaluation}
Let $U$ be a set with \textit{weight function}  $w\colon U\to \R^+$, and consider $\mathcal{P}_f (U)$, the set of finite subsets of $U$ partially ordered by inclusion. 

Then $w$ induces a length function on $(\mathcal{P}_f (U),\subseteq )$ by mapping a finite set $U'$ to $\sum_{x\in U'} w(x)$. This length function will also be denoted by $w$.
\end{Ex}
A length function  $l: (S,\leq ) \to (T,\lesssim)$  induces a function on finite chains 
 \[ l\colon Ch(S) \to Ch(T) ~~~~~~~~ X\mapsto l(X)=\{l(x)\mid x\in X\}. \] 
\begin{Def} \label{induced}
    Let $l\colon (S,\leq ) \to (T,\lesssim)$ be a length function.  
    
 The length function $l$  induces a \emph{chain length function}
 \[S\to Ch(T)\colon ~~~ x\longmapsto l^*(x):=max\{l(X)\mid X\in Ch(S);~ max(X)=x\}.\]

 If $X\in Ch(S)$ is a filtration and $l(X)=l^*(x)$, then $X$ is an \emph{$l^*$-filtration of x}.
 \end{Def} 
Note that by $(L2)$, the image of $l$ is totally ordered, so the lexicographical order on $l(Ch(S))$ is totally ordered too, so $l^*$ is well defined. In the above definition, it is enough to consider only filtrations of $x$ 
\begin{Ex}
    Let $U$ be a set. Consider $(\mathcal{P}_f (U), \subseteq)$ and the constant weight function $1$  sending each element to 1.
    
    Then for $U'\in \mathcal{P}_f (U)$ we have $1^*(U')=\{1,2,\dots, |U'|\}$.
    \end{Ex}
The following observation by Krause allows us to calculate $l^*$ recursively. 
\begin{Prop}[Krasue07, Prop. 1.5] \label{inductive}
  Let $l\colon (S,\leq ) \to (T,\lesssim)$ be a length function, $x\in S$.

    Then $l^*(x)=\{l(x)\}\cup max_{y<x} l^*(y)$
\end{Prop}
We end this section with the following theorem, investigating further the case from the example \ref{evaluation}. It is a crucial step in Theorem \ref{mainThm} about thin representations. Its proof is, however, purely combinatorial, so we present it in this more general form.
\begin{Thm} \label{filtration}
    Let $U$ be a set and consider a system of finite subsets $\mathcal{S}\subseteq  \mathcal{P}_f$ with a partial order  $\leq$ respecting the inclusion relation, i.e., if $X\leq X'$ then $X\subseteq X'$ for any $X, X'\in \mathcal{S}$. Let $X\in \mathcal{S}$ and  $\mathcal{X}= (X^1 < X^2 < \dots < X^m= X)$ be a filtration of $X$.

    Then there exists a weight function $w\colon U\to \R^+$ such that $\mathcal{X}$ is a $w^*$-filtration of $X$ for an induced length function $w\colon \mathcal{S}\to \R^+$.
\end{Thm}
Note that the following proof can be modified so that $\mathcal{X}$ is the only $w^*$-filtration of $X$. 
\begin{dukaz}
        We will use an induction on $m$. For each $n\leq m$, we construct a weight function  $w_n$ such that $\mathcal{X}^n$ is a $w_n^*$-filtration of $X^n$.  Then, $w_m$ will be the desired function. Because $\leq$ respects inclusions, the weights of the elements in $U\setminus X^m$ are irrelevant.

From the definition of filtration, the set $X^1$ is minimal in $(\mathcal{S},\leq)$, so it is its own $w^*$-filtratrion regardless of the weight function $w$. 

For $n<m$, assume that there is a weight function $w_n$ on $U$ such that $\mathcal{X}^n$ is a $w_n^*$-filtration of $X^n$.

We define $w_{n+1}$ as follows: Let $c$ be a real number such that $c> w(X^n)$. For each element $x\in X^{n+1}\setminus X^n$ we set $w_{n+1}(x)=c$. On the other elements, the weights stay the same. The chain $\mathcal{X}^n$ is thus a $w_{n+1}^*$-filtration of $X^n$ by the induction hypothesis. 

To prove that $\mathcal{X}^{n+1}$ is  a $w_{n+1}^*$-filtration of $X^{n+1}$, assume that
\[\mathcal{R}= (R^1 < R^2 <\dots < R^k = X^{n+1})\]
is a $w_{n+1}^*$-filtration of $X^{n+1}$. We will prove  $w(R^i)=w(X^i)$ for any $i\leq n+1$ and that $k=n+1$.

We work once again by induction; the case $i=1$ is clear from the definition. Now  assume that for any index $j< i$ we know $w_{n+1}(R^j)=w_{n+1}(X^j)$.  If  $R^i$ is contained in $X^n$, we can extend the chain $\mathcal{R}^i$  to a filtration $\mathcal{F}$ of $X^n$. If $w_{n+1}(R^i) < w_{n+1}(X^i)$ then we would have $w_{n+1}(\mathcal{F})>_L \mathcal{X}^n$. However,   $\mathcal{X}^n$ is  $w_{n+1}^*$-filtration of $X^n$, hence a contradiction. 

If   there is an element $x\in X^{n+1}\setminus X^n$ such that $x\in R^i$, then  $w_{n+1}(R^i)>c>w_{n+1}(X^n)$. From assumption  $w_{n+1}(\mathcal{R}^n)\geq_L w_{n+1}(\mathcal{X}^n)$ then follows that $i>n$ and $X^n$ is contained in $R^i$. By maximality of $X^n$ in $X^{n+1}$ we see that $R^i=R^k=X_{n+1}$.
\end{dukaz}

\section{Abelian length categories}\label{Abel}

This brief section recalls Krause's generalised definition of the Gabriel-Roiter measure from [12] and provides a minor reformulation that will be useful in the following sections.

\smallskip 

By a  \textit{length category}, we mean an essentially small abelian category all of whose objects are both artinian and noetherian, i.e., they admit a composition series.  By the Freyd-Mitchell embedding theorem,  such a category is essentially equivalent to a full subcategory of $Mod\text-R$ for some ring $R$.

The Gabriel-Roiter measures of an abelian category are then particular cases of induced chain length functions as defined in the previous section. But we limit ourselves only to length functions respecting short exact sequences in the following sense:

\begin{Def} \label{lfunction}
    Let $\mathcal{A}$ be a length category.  A \emph{length function on $\mathcal{A}$} is a map $l$ from the set of isomorphism classes of  objects of $\mathcal{A}$ to the set of nonnegative real numbers such that for any $X, Y, Z\in obj(\mathcal{A})$:
    
    (l1) $l([X])=0$ if and only if $X=0$.

(l2) If $0\rightarrow X \rightarrow Y \rightarrow Z \rightarrow 0$ is a short exact sequence, then $l([Y])=l([X])+l([Z])$.
\end{Def}

\begin{remark} \label{Jordan}
    By the Jordan-Hölder theorem, it is enough to specify the values of $l$ on the simple objects.
\end{remark}  

Let be $\mathcal{A}$ an abelian category and $l$ a length function on $\mathcal{A}$. Consider the poset of isomorphism classes of indecomposable objects, partially ordered by the subobject relation. A restriction of $l$ on this poset is a length function in the sense of Definition \ref{ldef}. Krause then defines the GR measure of an indecomposable object as the value of the induced chain length function on the corresponding isomorphism class. Using Proposition \ref{inductive} we reformulate this definition to a more computational-friendly version. 
\begin{Def}    
    Let  $X$ be an indecomposable object in a length category $\mathcal{A}$ and let $\mathcal{L}(X)$ be the set of its indecomposable subobjects, partially ordered by inclusion.
    
    Let $l$ be a length function on  $\mathcal{A}$. It induces a length function on $\mathcal{L}(X)$, also denoted $l$, via restriction. 
    
    The \emph{Gabriel-Roiter measure of $X$ with respect to $l$}, or the \emph{GR measure} for short,  is defined as $\mu_l(X):=l^*(X)$.

 An $l^*$-filtration of $X$ is called a \emph{GR filtration} with respect to $l$. An inclusion of two consecutive elements in a GR filtration is called a \emph{GR inclusion}.
\end{Def} 
Note that if $l$ is the composition length of an object, then this definition coincides with Gabriel's definition from [11].

\section{Quivers and thin representations} \label{ThinSec}

This section studies thin representations of quivers in terms of their support. Lemma \ref{thin} will allow us to work with a poset of quivers rather than in the category of representations. This is an integral step in the main theorem proved in the next section.  This transition to quivers also provides a simple combinatorial procedure for calculating the GR measure of a thin representation formulated in Algorithm \ref{Alg}.

\subsection{Quivers}

A \textit{quiver} $Q$ is a quadruple $(Q_0, Q_1,s,t)$ where $Q_0$ is a nonempty  set of \textit{vertices}, $Q_1$ a  set of \textit{arrows} and $s$ and $t$ are two maps $Q_1\to Q_0$ mapping an arrow to its \textit{source} and \textit{target}, respectively.  All quivers in this text are assumed to be finite, i.e., the sets $Q_0$ and $Q_1$ are finite. A quiver $Q'$ is a \textit{subquiver} of $Q$ if $Q'_0\subseteq Q_0$, $Q'_1\subseteq Q_1$ and the source  (target) map of $Q'$ is a restriction of a source (target) map of $Q$. In particular, each vertex corresponds to a \textit{trivial subquiver}. 

A multigraph obtained by forgetting the orientation of the arrows in $Q$ is the \textit{underlying graph} of $Q$. If the underlying graph of a quiver is a cycle on $n+1$ vertices ($n>0$), then the quiver is of type $\tilde{A}_n$. If it is a path on $n$ vertices, it is said to be of type $A_n$.

\begin{Def}
Let $Q$ be a connected quiver. 

For two connected subquivers  $Q'$ and $Q''$, we say that $Q'\leq_C Q''$ if $Q'$ is a subquiver of $Q''$ and if an arrow  $\alpha \in Q''_1$ has the source in $Q'_0$ then also the target of $\alpha$ is in $Q'_0$.

We denote by   $\mathcal{L}(Q)$ the set of connected subquivers of $Q$, partially ordered by $\leq_C$.
\end{Def}
Note that in this setting, theorem \ref{filtration} applies. We will only consider length functions given by weights of vertices.  
\begin{Def}
    Let $Q$ be a quiver. 
    
    A \emph{weight function} on $Q$ is a map $w\colon Q_0\to \R^+$. 
\end{Def}
 A weight function $w$ on the quiver $Q$ induces a length function on $\mathcal{L}(Q)$ mapping each connected subquiver to the sum of weights of its vertices. Proposition \ref{inductive} gives an algorithm to calculate $w^*(Q)$.
\begin{Alg}\label{Alg}
 First, we assign weights to all minimal elements in $\mathcal{L}(Q)$. If $Q$ is a tree, then these are only trivial quivers with one vertex - a sink in $Q$. It is enough to consider only those with minimal weight. From each minimal quiver, we start building a filtration. In each step, we chose a direct successor with the minimal weight. In cases where multiple successors have the same weight, we need to investigate all possibilities, as shown by the following example.
\end{Alg}
\begin{Ex}\label{Ex1}
Consider $Q= 1\leftarrow 2 \rightarrow 3  \leftarrow 4 \rightarrow 5 \leftarrow 6$ with a weight function $w$ prescribed as follows: $1\leftarrow 1 \rightarrow 1/2  \leftarrow 1 \rightarrow 1 \leftarrow 1$.

There are three minimal quivers (sinks), but vertex $3$ has the least weight, so any  $w^*$-filtration must start with $3$. There are two direct successors of $3$, both with weight $2.5$, namely $Q_l:= 1\leftarrow 2 \rightarrow 3 $ and $Q_r:= 3  \leftarrow 4 \rightarrow 5$. Starting with $Q_l$ we would get a filtration with weights $\{0.5,~2.5,~4.5,~5.5\}$ whereas starting with $Q_r$ we get a $w^*$-filtration with weights $\{0.5,~2.5,~3.5,~5.5\}$.

Thus, $w^*(Q)=\{0.5,~2.5,~3.5,~5.5\}$ with a $w^*$-filtration $(3<Q_r< 3  \leftarrow 4 \rightarrow 5 \leftarrow 6 <Q)$.
\end{Ex}

\subsection{Thin representations}

For a field $K$ and a quiver $Q$, a $K$-linear representation of $Q$ consists of a collection of finite-dimensional vector spaces $M_a$ for each $a\in Q_0$ and $K$-linear \textit{structural} maps $M_\alpha$ for each $\alpha\in Q_1$. The length category of all finite-dimensional $K$-representations of $Q$ is denoted by $rep_K(Q)$. A morphism $\Phi$ between two representations $N$ and $M$ is given by a collection of $K$-linear maps $\Phi_a\colon N_a\to M_a$ for each $a\in Q_0$ commuting with structural maps.

For $M\in rep_K(Q)$ indecomposable, we denote the poset of indecomposable subrepresentations by $\mathcal{L}(M)$. Following [19], we will call a filtration in $\mathcal{L}(M)$ an \textit{indecomposable filtration}.

  A $K$-representation of $Q$ is \textit{thin} if $dim_K(M_a)\leq 1$ for any $a\in Q_0$.  A thin representation such that all structural maps are identity maps will be called \textit{canonical}.  For $Q$ without oriented cycles, thin representations are exactly those representations that have finitely many subrepresentations. Bo Chen observed that if $M$ is a representation of a directed algebra (e.g., a representation of a  Dynkin quiver), then the first two elements of a GR filtration of $M$ are thin representations [5, Prop. 2.4.4]. If $N\hookrightarrow M$ is a GR inclusion of two indecomposable representations of a quiver of type $\tilde{A}_n$, then $M/N$  is a thin representation [7, Thm. 4.1].

\medskip 

The \textit{support of representation $M$} is a subquiver $supp(M)$ of $Q$ given by vertices and edges of $Q$ such that the corresponding linear spaces and structural maps are nonzero. Note that a structure map $\phi_\alpha$ in thin representation is either the zero map ($\alpha \notin supp(M)_1)$ or an isomorphism ($\alpha \in supp(M)_1$).

The following lemma will allow us to work with supports instead of thin representations.  The presented results are probably not novel. The case when the quiver is of type $A_n$ is well-known in the persistence theory; see [13]. However, the more general case needed here is not covered in the literature. Thus, we outline the proofs for the reader's convenience.  Some results are formulated more generally, as they can be obtained with little or no changes in proof.

\begin{Lemma} \label{thin}
    Let $Q$ be an acyclic quiver, let $N, N', M$ be thin representations of $Q$ and suppose that the underlying graph of $supp(N)$ contains no cycles.  

(1) $M$ is indecomposable iff $supp(M)$ is connected.

(2) $N\cong N'$ iff $supp(N)=supp(N')$.

(3) Suppose that $N$ and $M$ are indecomposable. Then there is a monomorphism from $N$ to $M$ iff $supp(N)\leq_C supp(M)$.     
\end{Lemma} The forward implications in (1) and (2) hold in general. 
\begin{dukaz}
To prove (1), observe that if a thin representation is decomposable, then its decomposition into a direct sum of indecomposable summands induces a decomposition of the quiver $supp(M)$.

From (1), it follows that it is enough to prove (2) for the case when $supp(N)$ is connected. In this case, the underlying graph of $supp(N)$ is a tree. 

Any thin representation $N$ whose support is a tree is isomorphic to a canonical representation. We work by induction on the number of vertices. For a leaf $a$ in $supp(N)$, we can create a new representation $N^a$  from $N$ by setting $N_a=0$ and $N_\alpha =0$ for the arrow whose source or target is $a$. By induction, $N^a$ is isomorphic to a canonical representation.  Because there is only one arrow adjacent to $a$, we can extend this isomorphism to an isomorphism defined on the whole $N$. Thus, both $N$ and $N'$ are isomorphic to the canonical representation of $supp(N)=supp(N')$.

For the forward implication in (3), let $\phi:N\to M$ be a monomorphism and $\alpha\in supp(M)_1$ such that $N_{t(\alpha)}=0$. Because components of $\phi$ commute with structrual maps we see that $M_\alpha\circ \phi_{s(\alpha)}=0$, i.e., $Im(\phi_{s(\alpha)})\subseteq Ker(M_\alpha)$. But for thin representation, the assumption $M_{t(\alpha)}\neq 0$  implies that  $Ker(M_\alpha)=0$, hence $\phi_{s(\alpha)}$ is the zero map. But $\phi$ is a monomorphism, hence $N_{s(\alpha)}=0$. So any arrow with the source in $supp(N)_0$ also has the target in $supp(N)_0$.

For the other implication, consider a subrepresentation $M'$ of $M$  such that $supp(M')=supp(N)$ and all the nonzero structure maps from $M'$ coincide with structure maps from $M$.
Because $supp(M')$ is a tree, using (2), we get an isomorphism $N\cong M'$. Composing it with the inclusion $M'\leq M$, we get a monomorphism from $N$ to  $M$.
\end{dukaz}

\section{Generalised GR measure and thin representations} \label{mainSec}

This section shows that the GR measure of a thin representation of a Dynkin or Euclidean quiver can be calculated as an induced chain length function of its support quiver. Proposition \ref{inductive} thus gives an easy combinatorial procedure for calculating GR measures. This also shows that GR measures and the corresponding GR filtrations of thin representations are field-independent. To illustrate this approach, a result by Bo Chen is generalised in Example \ref{injectivefactor}. As a corollary of theorem \ref{filtration}, it follows that any indecomposable filtration is a GR filtration for some length function. The section ends with an example showing the limits of presented results.

\medskip

Let $Q$ be an \textit{acyclic quiver}, i.e., a quiver containing no oriented cycles. For a field $K$, a \textit{path algebra KQ} is the $K$-vector space with a basis consisting of paths in $Q$. The product of two paths is a  concatenation when possible and zero otherwise. Path algebras are hereditary, and if $K$ is algebraically closed, any hereditary finite-dimensional $K$-algebra is Morita equivalent to $KQ$ for some acyclic quiver $Q$. See [2] for details. 

The category of finitely generated right $KQ$-modules is denoted by $mod\text- KQ$. Algebra $KQ$ is finite-dimensional, and hence artinian as a ring. The category $mod\text- KQ$ is a length category equivalent to $rep_K(Q)$. It has only finitely many simple representations (up to isomorphism) corresponding to a thin representation $\{S(i)\mid i\in Q_0\}$ whose supports are trivial subquivers.

If $l$ is a length function on  $mod\text- KQ$, it induces a weight function on $Q$ by the formula $w(i)=l(S(i))$ for any $i\in Q_0$.  Then, for any representation $M\in mod\text- KQ$, we have $l(M)=\sum_{i\in Q_0} dim_K(M_i)\cdot l(S(i))$. In particular, the length of a thin representation corresponds to the weight of its support. 

\begin{Prop}\label{translate}
Let $Q$ be acyclic quiver, $K$ a field and $M\in rep_K(Q)$ a thin representation such that the underlying graph of $supp(M)$ is either a  tree or a cycle. Let $l$ be a length function on $mod\text- KQ$ and $w$ its induced weight function on  $\mathcal{L}(Q)$. 

    Then  $\mu_l(M)=w^*(supp(M))$ and if $\mathcal{Q}$ is a $w^*$-filtration of $supp(M)$, then the canonical representations of the quivers from $\mathcal{Q}$ form a GR filtration of $M$. 
\end{Prop}
In particular, if $Q$ is a Dynkin or Euclidean quiver, the proposition applies to all thin representations. 
\begin{dukaz}
    If the underlying graph of $supp(M)$ is a tree, by lemma \ref{thin} the map \textit{supp} assigning to each representation its support restricts to a poset isomorphism between the poset of isomorphism classes of subrepresentations of $M$ and $\mathcal{L}(supp(M))$, such that $l(N)=w(supp(N))$ for any $N\leq M$. 
    
     When $supp(M)$ is of type $\tilde{A}_n$, observe that the underlying graph of a support of any nontrivial subrepresentation of $M$  is a tree. Thus if $N\hookrightarrow M$ is a GR-inclusion, then $\mu_l(M)=\{l(M)\}\cup \mu_l(N)=\{l(M)\} \cup w^*(supp(N))$ and $supp(M)$ is a successor of $supp(N)$ in $\mathcal{L}(Q)$.
\end{dukaz}
This proposition allows us to work only with weighted quivers instead of thin representations. Algorithm \ref{Alg} then gives us a simple method to calculate the GR measure of any thin representation. In particular, for quivers of type $A_n$ (path on $n$-vertices), all Gabriel-Roiter measures can be calculated this way.
\begin{Ex}\label{Ex2}
    Consider quiver $Q$ from Example \ref{Ex1} and a $K$-representation \[M=0\xleftarrow{\cdot 0} 0 \xrightarrow{\cdot 0} K  \xleftarrow{\cdot 3} K \xrightarrow{id} K \xleftarrow{\cdot 2} K.\] A $K$-dimension of a representation is a length function on $mod\text-KQ$.

    Then $supp(M)= 3  \leftarrow 4 \rightarrow 5\leftarrow 6$ and the weight function $w$ induced by a $K$-dimension assigns 1 to each vertex. We now follow Algorithm $\ref{Alg}$, to calculate $\mu(M)$ and the corresponding GR-filtrations.

There are two minimal quivers (sinks), both with weights 1. There is only one possible filtration in $\mathcal{L}(supp(M))$ starting with 3, that is $(3<\ 3  \leftarrow 4 \rightarrow 5<supp(M))$.

Two options start with 5, but they have different weights. This gives a filtration $(5< 5\leftarrow 6 <supp(M))$.

We conclude that $\mu(M)=\{1,2,4\}$ and there is only one GR-filtration of $M$, namely 
\[ (S(5) <0\xleftarrow{\cdot 0} 0 \xrightarrow{\cdot 0} 0  \xleftarrow{\cdot 0} 0 \xrightarrow{\cdot 0} K \xleftarrow{\cdot 2} K < M).\]
\end{Ex}

The following example shows a further application of this method. We generalise a known result by Bo-Chen, [9, Prop. 3.2], which assumes $K$ to be algebraically closed.
\begin{Ex}\label{injectivefactor}
Let $Q$ be a quiver of type $\tilde{A}_n$ and let $M$ be a quasi-simple homogenous representation. This is equivalent to $supp(M)=Q$. 

Let  $N\leq M$
 be  a maximal indecomposable submodule. Then  $supp(N)$ is maximal in $supp(M)$ with respect to $\leq_C$. Then $supp(M)\setminus supp(N)$ consists of one vertex $i$ and two arrows starting in $i$. We see that $M/N\cong S(i)$. Because $i$ is a source, the module $S(i)$ is injective; see [1, Chapter III.]. 
\end{Ex} 
 \begin{Thm}\label{mainThm} Let $K$ be a field, $Q$ an acyclic quiver.   Let $M\in rep_K(Q)$ be a thin indecomposable representation such that $supp(M)$ is either a tree or $\tilde{A}_n$. Let $\mathcal{M}=(M_1 < M_2 < \dots < M_n = M)$ be an indecomposable filtration of $M$.

 Then there exists a length function $l$ on $mod\text- KQ$ such that $\mathcal{M}$ is a GR-filtration of $M$ with respect to $l$.     
 \end{Thm}
 \begin{dukaz}
     By Proposition \ref{translate}, it is enough to prove that there exists a weight function $w$ on $supp(M)$ such that $(supp(M_1) <_C supp(M_2) <_C \dots <_C supp(M_n) = supp(M))$ is a $w^*$-filtration. But that follows directly from Theorem \ref{filtration}.
 \end{dukaz}
\begin{Ex}
    Recall quiver $Q$ and representation $M$ as in Example \ref{Ex2} and consider an indecomposable filtration 
    \[\mathcal{M}=(S(3) < 0\xleftarrow{\cdot 0} 0 \xrightarrow{\cdot 0} K  \xleftarrow{\cdot 3} K \xrightarrow{id} K \xleftarrow{\cdot 0} 0<M).\]

    We aim to find a length function $l$ on $mod\text-KQ$ such that $\mathcal{M}$ is a GR-filtration with respect to $l$. By Proposition \ref{translate}, this amounts to finding  a weight function $w$  on $Q$ such that $(3<3  \leftarrow 4 \rightarrow 5<supp(M))$ is a $w^*$-filtration.

Following the proof of Theorem \ref{filtration}, we find a sequence of weight functions $w_1, w_2$ and $w_3$, such that $\mathcal{M}^i$ is a $w_i^*$-filtration for $i\in\{1,2,3\}$ and set $w:=w_3$. Let $w_1$ be constantly 1 and 
\begin{gather*}
    w_2:= 1\leftarrow 1 \rightarrow 1  \leftarrow 2 \rightarrow 2 \leftarrow 1\\
    w_3:= 1\leftarrow 1 \rightarrow 1  \leftarrow 2 \rightarrow 2 \leftarrow 6.
\end{gather*}

Translating back to $mod\text- KQ$, we get $l(S(i))=1$ for $i\in\{1,2,3\}$, $l(S(4))=l(S(5))=2$ and $l(S(6))=6$. By Remark \ref{Jordan}, this completely determines a length function on $mod\text- KQ$. Filtration $\mathcal{M}$ is then a (unique) GR-filtration of $M$ with respect to $l$.

Note that choosing the weight function to be $w':=  1\leftarrow 1 \rightarrow 1 \leftarrow 2 \rightarrow 1 \leftarrow 6$, we would get two GR-filtrations of $M$, namaely $\mathcal{M}$ and  
\[(S(5) < 0\xleftarrow{\cdot 0} 0 \xrightarrow{\cdot 0} K  \xleftarrow{\cdot 3} K \xrightarrow{id} K \xleftarrow{\cdot 0} 0<M).\]    
\end{Ex}
 
 The following example shows the limits of the theorem when the representations are not thin.  
\begin{Ex}\label{D4}
    Let $\xymatrix@=10pt{
Q:  & 1 \\
2 & 3 \ar[l] \ar[u] \ar[r] & 4   
    }   
   $ be a quiver of type $D_4$. It is a quiver of the Dynkin type, so all representations (up to isomorphism) are determined by their dimension vectors. The Auslander-Reiten quiver is
    
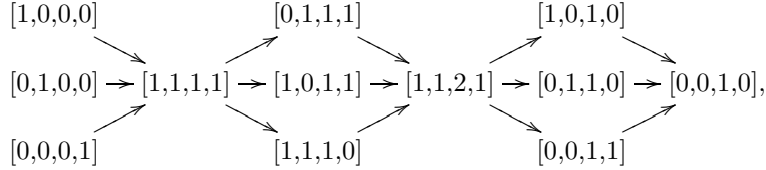
\[\xymatrix@=10pt{
[1,0,0,0] \ar[rd]  & &  [0,1,1,1] \ar[rd] && [1,0,1,0] \ar[rd]\\
[0,1,0,0]   \ar[r] & [1,1,1,1] \ar[ru] \ar[r] \ar[rd]&[1,0,1,1] \ar[r]& [1,1,2,1] \ar[ru] \ar[r] \ar[rd] &[0,1,1,0] \ar[r] & [0,0,1,0],\\
[0,0,0,1] \ar[ru] && [1,1,1,0] \ar[ru]&&[0,0,1,1] \ar[ru]
}\]
where in each vertex, there is a dimension vector of indecomposable representations in the corresponding equivalence class.

    Let $M$ be an indecomposable representation with $dim(M)=[1,1,2,1]$ and let $N$ be an indecomposable representation of composition length four. Then $\mathcal{M}_i = (S(i)\leq N \leq M)$ is an indecomposable filtration of $M$ for $i\in \{1,2,4\}$. We show that $\mathcal{M}_i$ cannot be a GR filtration.

    Let $j$ be a sink such that $l(S(j))\leq l(S(i))$ for any $i\in \{1,2,4\}$. There are (up to isomorphism) three indecomposable representations of length three. They are all thin. We choose among them $N'$ such that $j\in supp(N)_0$.
    
    Representation $N'$ injects into $M$, but this morphism cannot be factored through $N$; hence $(S(j)\leq N'\leq M)$ is also an indecomposable filtration.
    
    Regardless of the choice of the length function $l$, we always have $l(N')<l(N)$.  Hence 
    \[\mu_l(N)=\{l(S(j)), l(N)\} < \{l(S(j)), l(N')\} = \mu_l(N')\]
\end{Ex}
\begin{remark}
    If $K$ is a finite field, Theorem \ref{filtration} gives us a length function such that $(S(i)\leq N \leq M)$ is a $l^*$-function. Thus, such a length function does not satisfy Definition \ref{lfunction}.
\end{remark}

\noindent\textbf{Aknowledgements:}
 The author wishes to express his thanks and appreciation to Doc. Jan Šťovíček for introducing him to the concept of the GR measure and for his guidance and valuable suggestions during the research. 

\medskip 

\noindent \textbf{Author info:}

RNDr. Dominik Krasula

\noindent Address:

Czech Republic

Brantice 281

793 93 Brantice

\noindent e-mail: krasula@karlin.mff.cz

\medskip 

\noindent \textbf{Bibliography}

 [1] Assem, I. Skowronski, A.  Simson, D. (2006). Elements of the Representation Theory of Associative Algebras: Techniques of Representation Theory.  \textit{LMS Student Texts}. Cambridge: Cambridge University Press. vol. 65. doi:10.1017/CBO9780511614309

\smallskip 

 [2] Auslander, M. Reiten, I. Smal\o, S. (1995). Representation Theory of Artin Algebras. \textit{Cambridge Studies in Advanced Mathematics}. Cambridge University Press.  vol. 36. doi:10.1017/CBO9780511623608

\smallskip

 [3] Auslander, I. Smal\o, S. (1980). Preprojective Modules over Artin Algebras. \textit{J. Algebra}.  66(1). 61-122. https://doi.org/10.1016/0021-8693(80)90113-1

\smallskip

 [4] Bahlekeh, A., Fotouhi, F.S.,  Salarian, S. (2018). Representation-theoretic properties of balanced big Cohen–Macaulay modules. \textit{Mathematische Zeitschrift}, 293, 1673-1709. https://doi.org/10.1007/s00209-019-02257-1

\smallskip

 [5] Chen, B. (2006). The Gabriel-Roiter Measure For Representation-Finite Hereditary Algebras. [Doctoral dissertation, Der Universität Bielefeld].

\smallskip

  [6] Chen, B. (2007). The Gabriel–Roiter measure for representation-finite hereditary algebras. \textit{J. Algebra}. 309(1). 292-317. https://doi.org/10.1016/j.jalgebra.2006.08.018

\smallskip

 [7] Chen. B. (2007). The Auslander-Reiten sequences ending at Gabriel-Roiter factor modules over tame hereditary algebras. \textit{J. Algebra its Appl.} 6(6). 951-963. https://doi.org/10.1142/S0219498807002570

\smallskip

 [8] Chen, B. (2008) Comparison of Auslander–Reiten Theory and Gabriel–Roiter Measure Approach to the Module Categories of Tame Hereditary Algebras.  \textit{Commun. Algebra}. 36(11). 4186-4200.  https://doi.org/\allowbreak10.1080/00927870802175113

\smallskip

 [9] Chen, B. (2008). The Gabriel–Roiter measure for $\tilde{A}_n$. \textit{J. Algebra}. 320(7). 2891-2906. https://doi.org/\allowbreak10.1016/j.jalgebra.2008.06.022

\smallskip

 [10] Chen, B. (2010). The Gabriel-Roiter submodules of simple homogeneous modules. \textit{Proc. Amer. Math. Soc.} 138. 3415-3424. https://doi.org/10.1090/S0002-9939-2010-10243-5

\smallskip

 [11] Gabriel, P. (1973). Indecomposable representations II. \textit{Symposia Mathematica}. 11. 81-104.

\smallskip

[12] Krause, H. (2007). An axiomatic characterization of the Gabriel-Roiter measure.  \textit{Bull. London Math.} 39(4). 550-558. https://doi.org/10.1112/blms/bdm033

\smallskip

 [13] Oudot, S. (2015). Persistence Theory - From Quiver Representations to Data Analysis. \textit{}\textit{Mathematical Surveys and Monographs}. Amer. Math. Soc., Providence.  vol. 209. https://doi.org/10.1090/surv/209

\smallskip

[14] Ringel, C. M. (2005). The Gabriel–Roiter measure. \textit{Bulletin des Sciences Mathématiques.} 129(9). 726-748. https://doi.org/10.1016/j.bulsci.2005.04.002

\smallskip 

 [15] Ringel, C.M. (2006). The theorem of Bo Chen and Hall polynomials  \textit{Nagoya Math. J}. 183. 143-160. https://doi.org/10.1017/S0027763000009284

\smallskip

 [16] Ringel, C. M. (2006). Foundation of the representation theory of Artin algebras, Using the Gabriel-Roiter measure. In:\textit{ Trends in representation theory of algebras and related topics}. Contemporary Mathematics. Amer. Math. Soc., Providence, vol. 406. 105-135.  https://doi.org/10.1090/conm/406/07656. 

\smallskip

 [17] Ringel, C. M. (2010). Gabriel–Roiter inclusions and Auslander–Reiten theory.  \textit{J. Algebra.
} 324(12). 3579-3590. https://doi.org/10.1016/j.jalgebra.2010.09.003

\smallskip

 [18] Roiter, A. V. (1968). Unbounded dimensionality of indecomposable representations
of an algebra with an infinite number of indecomposable representations. \textit{Mathematics of the USSR-Izvestiya}. 2(6). 1223-1230. https://doi.org/10.1070/IM1968v002n06ABEH000727

\smallskip

 [19] Schmidmeier, M.  Tyler, H. R. (2014). The Auslander–Reiten
Components in the Rhombic Picture. \textit{Commun. Algebra}. 42(3). 1312-1336. https://doi.org/10.1080/00927872.2012.738341

\smallskip

 [20] Szántó, C.,   Szöllősi, I. (2014). Hall polynomials and the Gabriel–Roiter submodules of simple homogeneous modules. \textit{Bull. London Math}. 47(2). 206-216. https://doi.org/10.1112/blms/bdu109

\end{document}